\documentclass[11pt]{article}
\usepackage{latexsym,amssymb,amsmath,geometry,cite,enumerate,float,verbatim}
\newtheorem{theorem}{Theorem}[section]
\newtheorem{lemma}[theorem]{Lemma}
\newtheorem{conjecture}[theorem]{Conjecture}
\newtheorem{corollary}[theorem]{Corollary}

\usepackage{lineno}
\usepackage{setspace}

\begin{document}

\onehalfspace

\title{Equating $k$ Maximum Degrees in Graphs without Short Cycles}

\author{M. F\"{u}rst$^1$\and 
M. Gentner$^1$\and 
M.A. Henning$^2$\and 
S. J\"{a}ger$^1$\and 
D. Rautenbach$^1$}

\date{}

\maketitle

\begin{center}
{\small 
$^1$ Institute of Optimization and Operations Research, Ulm University, Germany\\
\texttt{maximilian.fuerst,michael.gentner,simon.jaeger,dieter.rautenbach@uni-ulm.de}\\[3mm]
$^2$ Department of Pure and Applied Mathematics, University of Johannesburg, South Africa\\
\texttt{mahenning@uj.ac.za}
}
\end{center}

\begin{abstract}
For an integer $k$ at least $2$, and a graph $G$,
let $f_k(G)$ be the minimum cardinality of a set $X$ of vertices of $G$
such that $G-X$ has either $k$ vertices of maximum degree or order less than $k$.
Caro and Yuster (Discrete Mathematics 310 (2010) 742-747) conjectured that, 
for every $k$, 
there is a constant $c_k$ such that $f_k(G)\leq c_k \sqrt{n(G)}$ for every graph $G$.
Verifying a conjecture of Caro, Lauri, and Zarb (arXiv:1704.08472v1), 
we show the best possible result that, 
if $t$ is a positive integer, 
and $F$ is a forest of order at most $\frac{1}{6}\left(t^3+6t^2+17t+12\right)$,
then $f_2(F)\leq t$.
We study $f_3(F)$ for forests $F$ in more detail obtaining similar almost tight results,
and we establish upper bounds on $f_k(G)$ for graphs $G$ of girth at least $5$.
For graphs $G$ of girth more than $2p$, for $p$ at least $3$, 
our results imply $f_k(G)=O\left(n(G)^{\frac{p+1}{3p}}\right)$.
Finally, we show that, for every fixed $k$, 
and every given forest $F$, 
the value of $f_k(F)$ can be determined in polynomial time.
\end{abstract}

{\small 

\begin{tabular}{lp{13cm}}
\end{tabular}
}

\pagebreak

\section{Introduction}

Every finite, simple, and undirected graph has at least two vertices of equal degree,
and this lower bound on the number of repeated degrees can be improved for restricted graph classes \cite{cw}.
Caro, Shapira, and Yuster \cite{csy} proved the surprising result that, 
for every positive integer $k$,
there is a constant $c_k$ such that,
for every graph $G$, there is a set $X$ of at most $c_k$ vertices such that 
$G-X$ has at least $\min\{ k,n(G)-|X|\}$ many vertices of equal degree,
where $n(G)$ denotes the order of $G$.

In \cite{cy} Caro and Yuster considered an analogous problem for the maximum degree.
For an integer $k$ at least $2$, and a graph $G$,
let $f_k(G)$ be the minimum cardinality of a set $X$ of vertices of $G$
such that $G-X$ has either $k$ vertices of maximum degree or order less than $k$.

Caro and Yuster pose the following intriguing conjecture.

\begin{conjecture}[Caro and Yuster \cite{cy}]\label{conjecturecy}
For every integer $k$ at least $2$, there is a constant $c_k$ such that
$f_k(G)\leq c_k \sqrt{n(G)}$ for every graph $G$.
\end{conjecture}
They describe graphs $G$ with $f_2(G)\geq (1-o(1))\sqrt{n(G)}$
showing that the upper bound in Conjecture \ref{conjecturecy} has the best possible growth rate,
that is, forcing many vertices of maximum degree is considerably harder than forcing many vertices of equal degree.
Furthermore, they verify the conjecture for $k\in \{ 2,3\}$
by showing that $c_2=\sqrt{8}$ and $c_3=43$ have the desired properties.
They also prove the following result, which implies the conjecture for $C_4$-free graphs.

\begin{theorem}[Caro and Yuster \cite{cy}]\label{theoremcy}
Let $k$ and $t$ be positive integers at least $2$.
If $G$ is a $K_{2,t}$-free graph of order at least $t^2{k\choose 2}^2$, 
then $f_k(F)\leq (3k-3)\sqrt{n(G)}$.
\end{theorem}
In \cite{clz} Caro, Lauri, and Zarb show that $\sqrt{2}$ is the best possible value for $c_2$,
and, for forests $F$, they improve the growth rate of the upper bound on $f_k(F)$ from the second to the third root of the order as follows.

\begin{theorem}[Caro, Lauri, and Zarb \cite{clz}]\label{theoremclz}
If $k$ is an integer at least $2$, and $F$ is a forest of order at least $(2k-1)^3$,
then $f_k(F)\leq (2k-1)n(F)^{\frac{1}{3}}$.
\end{theorem}
For $k=2$, they formulate a precise conjecture,
and construct graphs showing that their conjecture would be tight.

\begin{conjecture}[Caro, Lauri, and Zarb \cite{clz}]\label{conjecture1}
If $t$ is a positive integer, and $F$ is a forest of order at most $\frac{1}{6}\left(t^3+6t^2+17t+12\right)$,
then $f_2(F)\leq t$.
\end{conjecture}
In the present paper we show this conjecture.
Furthermore, we study $f_3(F)$ for forests $F$ in more detail obtaining almost tight results,
and we give improved upper bounds on $f_k(G)$ for graphs $G$ of girth at least $5$.
For graphs $G$ of girth more than $2p$, for $p$ at least $3$, 
our results imply $f_k(G)=O\left(n(G)^{\frac{p+1}{3p}}\right)$,
and we obtain considerable improvements of Theorem \ref{theoremclz}.
Finally, we show that, for every fixed integer $k$ at least $2$, 
and every given forest $F$, 
the value of $f_k(F)$ can be determined in polynomial time.

The influence of degree multiplicities on graph parameters or
large sets of vertices of equal degree satisfying additional properties
have been studied in several papers such as \cite{ab,bs}; 
see \cite{clz} for further discussion.
Before we proceed to our results, we collect some notation.
Let $G$ be a graph.
The size of $G$ is denoted by $m(G)$.
For a vertex $u$ of $G$, the degree of $u$ is denoted by $d_G(u)$.
The maximum degree of $G$ is denoted by $\Delta(G)$.

For an integer $n$, let $[n]$ be the set of all positive integers at most $n$.
If $G$ has order $n$ and degree sequence $d_1\geq d_2 \geq \ldots \geq d_n$, 
then let $\Delta_i(G)$ be $d_i$ for $i\in [n]$; in particular, 
$\Delta_1(G)$ is the maximum degree of $G$,
and $G$ has at least $k$ vertices of maximum degree if and only if $\Delta_1(G)=\Delta_k(G)$.

\section{Upper bounds}

Our first goal is the proof of Conjecture \ref{conjecture1}.
The following result from \cite{clz} was the key insight needed to obtain the best possible value for $c_2$.

\begin{theorem}[Caro, Lauri, and Zarb \cite{clz}]\label{lemma1}
If $t$ is a positive integer, and $G$ is a graph with $\Delta(G)\leq {t+2\choose 2}$, then $f_2(G)\leq t$.
\end{theorem}
Since a forest has less edges than vertices, 
the following result immediately implies Conjecture \ref{conjecture1}.

\begin{theorem}\label{theorem1}
If $t$ is a positive integer, and $F$ is a forest of size less than $\frac{1}{6}\left(t^3+6t^2+17t+12\right)$, 
then $f_2(F)\leq t$.
\end{theorem}
{\it Proof:} For a positive integer $t$, let $n(t)=\frac{1}{6}\left(t^3+6t^2+17t+12\right)$.
The proof is by induction on $t$.
Let $\Delta_i=\Delta_i(F)$ and let $u_i$ be such that $d_F(u_i)=\Delta_i$ for $i\in [2]$,
where $u_1$ and $u_2$ are distinct.

For $t=1$, we have $m(F)\leq n(1)-1=5$.
Clearly, we may assume that $f_2(F)>0$, that is, $\Delta_1>\Delta_2\geq 1$.
If $\Delta_2=1$, then $F$ is the union of a star $K_{1,\Delta_1}$ and copies of $K_1$ and $K_2$,
and removing $u_1$ yields two vertices of maximum degree $0$ or $1$.
Hence, we may assume that $\Delta_2\geq 2$, 
which, using $m(F)\leq 5$, implies that $3\leq \Delta_1\leq 4$.
If $\Delta_1=3$, then removing a neighbor of $u_1$ that does not lie in $N_F[u_2]$ 
yields two vertices of maximum degree $2$.
Note that such a neighbor exists, because $F$ is a forest.
Hence, we may assume that $\Delta_1=4$,
which implies that $F$ arises by subdividing one edge of a star $K_{1,4}$ once,
and removing $u_1$ yields two vertices of maximum degree $1$.

Now, let $t\geq 2$.
If $\Delta_1\leq {t+2\choose 2}$, then Theorem \ref{lemma1} implies $f_2(F)\leq t$.
Hence, we may assume that $\Delta_1\geq {t+2\choose 2}+1$.
If $F'=F-u_1$, then 
\begin{eqnarray*}
m(F') & = & m(F)-\Delta_1\\
& < & \frac{1}{6}\Big(t^3+6t^2+17t+12\Big)-\Big(\frac{1}{2}t^2+\frac{3}{2}t+2\Big)\\
& = & \frac{1}{6}\Big((t-1)^3+6(t-1)^2+17(t-1)+12\Big)\\
& = & n(t-1).
\end{eqnarray*}
By induction, we obtain $f_2(F)\leq 1+f_2(F')\leq 1+(t-1)=t$,
which completes the proof. $\Box$

\medskip

\noindent In order to better understand $f_k(F)$ for forests $F$, we first consider the case $k=3$.

Our next result suitably generalizes Theorem \ref{lemma1}. 

\begin{theorem}\label{theorem2}
If $t$ is an integer at least $2$, and $F$ is a forest with $\Delta_1(F)+2\Delta_2(F)\leq {t+2\choose 2}+2$,
then $f_3(F)\leq t$.
\end{theorem}
{\it Proof:} The proof is by induction on $t$.
Clearly, we may assume that $F$ has at least three vertices,
and that $\Delta_1(F)>\Delta_3(F)$.
Let $\Delta_i=\Delta_i(F)$ and let $u_i$ be such that $d_F(u_i)=\Delta_i$ for $i\in [3]$,
where $u_1$, $u_2$, and $u_3$ are distinct.

For $t=2$, we have $\Delta_1+2\Delta_2\leq {2+2\choose 2}+2=8$.
If $\Delta_1=1$, then $F$ is the union of copies of $K_1$ and $K_2$,
and removing one vertex of degree $1$
yields either three vertices of maximum degree $0$,
or a graph with less than $3$ vertices.
Hence, we may assume that $\Delta_1\geq 2$.
If $\Delta_2=1$, then $F$ is the union of a star $K_{1,\Delta_1}$,
copies of $K_1$, and $p$ copies of $K_2$.
If $p=0$ or $p\geq 2$, then let $X=\{ u_1\}$, and,
if $p=1$, then let $X$ contain $u_1$ and exactly one vertex from the unique $K_2$ component.
It is easy to check that $F-X$ has three vertices of maximum degree.
Hence, we may assume that $\Delta_2\geq 2$, 
which, using the upper bound on $\Delta_1+2\Delta_2$, 
implies $\Delta_1\in \{ 2,3,4\}$ and $\Delta_2=2$.
First, we assume that $\Delta_1=2$.
Clearly, we may assume that $\Delta_3=1$.
If $u_1$ and $u_2$ are non-adjacent, then $F$ contains two copies of $P_3$,
and removing one endvertex from each copy yields four vertices of maximum degree $1$.
If $u_1$ and $u_2$ are adjacent, let $F$ contain $p$ $K_2$ components.
If $p=0$, then let $X=\{ u_1,u_2\}$, and,
if $p\geq 1$, then let $X=\{ u_1\}$.
It is easy to check that either $n(F-X)<3$ or 
$F-X$ has three vertices of maximum degree.
Hence, we may assume that $\Delta_1\geq 3$.
If $\Delta_3=2$, then removing $\Delta_1-2$ neighbors of $u_1$
that do not belong to $N_F[u_2]\cup N_F[u_3]$ 
yields three vertices of maximum degree $2$.
Hence, we may assume that $\Delta_3=1$.
If $F$ has no $K_2$ component, then removing $u_1$ and $u_2$ yields
three vertices of maximum degree $0$.
Hence, we may assume that $F$ has a $K_2$ component.
If $u_1$ is adjacent to $u_2$, then let $X=\{ u_1\}$, and, 
if $u_1$ is non-adjacent to $u_2$, then let $X$ contain $u_1$ and exactly one neighbor of $u_2$.
It is easy to check that 
$F-X$ has three vertices of maximum degree.

Now, let $t\geq 3$.
First, suppose that $\Delta_1+\Delta_2-2\Delta_3\leq t$.
Clearly, we may assume that $\Delta_3\geq 1$.
If $\Delta_3=1$, then 
either removing $u_1$ and $u_2$
or removing $\Delta_1-1$ neighbors of $u_1$ that do not belong to $N_F[u_2]$
and $\Delta_2-1$ neighbors of $u_2$ that do not belong to $N_F[u_1]$
yields three vertices of maximum degree $0$ or $1$.
Hence, we may assume that $\Delta_3\geq 2$.
Now, removing 
$\Delta_1-\Delta_3$ neighbors of $u_1$ that do not belong to $N_F[u_2]\cup N_F[u_3]$
and 
$\Delta_2-\Delta_3$ neighbors of $u_2$ that do not belong to $N_F[u_1]\cup N_F[u_3]$
yields three vertices of maximum degree $\Delta_3$.
Again, all these vertices exist, because $F$ is a forest.
Hence, we may assume that $\Delta_1+\Delta_2-2\Delta_3\geq t+1$.

Let $F'=F-u_1$. 
Clearly, $\Delta_1(F')\leq \Delta_2$ and $\Delta_2(F')\leq \Delta_3$.
If $\Delta_2+2\Delta_3\leq {t+1\choose 2}+2$,
then $\Delta_1(F')+2\Delta_2(F')\leq {t+1\choose 2}+2$, and, by induction, 
$f_3(F)\leq 1+f_3(F')\leq 1+(t-1)=t$.
Hence, we may assume that $\Delta_2+2\Delta_3\geq {t+1\choose 2}+3$,
and, we obtain 
\begin{eqnarray*}
\Delta_1+2\Delta_2 & = & \Big(\Delta_1+\Delta_2-2\Delta_3\Big)+\Big(\Delta_2+2\Delta_3\Big)\\
& \geq & \Big(t+1\Big)+\left({t+1\choose 2}+3\right)\\
& = & {t+2\choose 2}+3,
\end{eqnarray*}
which is a contradiction.
$\Box$

\medskip

\noindent Since $f_3(K_{1,5}\cup P_3\cup P_3)=3$, 
the base case of the induction in the previous proof is best possible.

Note that $f_3(K_{1,3}\cup K_2)=2$ 
shows that Theorem \ref{theorem2} is not true for $t=1$.

By a simple inductive argument,
Theorem \ref{theorem2} implies a lower bound on the sum of the largest degrees in terms of $f_3(F)$.

\begin{corollary}\label{corollary1}
If $t$ is an integer at least $2$, and $F$ is a forest with $f_3(F)>t$, then 
\begin{enumerate}[(i)]
\item $\Delta_t(F)\geq 2$, 
\item $\Delta_{t+1-i}(F)+2\Delta_{t+2-i}(F)\geq {i+2\choose 2}+3$ for every $i\in [t]\setminus \{ 1\}$, and
\item $\Delta_1(F)+\Delta_2(F)+\cdots+\Delta_t(F)\geq \frac{1}{18}t^3 + \frac{1}{3}t^2 + \frac{29}{18}t$.
\end{enumerate}
\end{corollary}
{\it Proof:} Let $\Delta_i=\Delta_i(F)$ and let $d_F(u_i)=\Delta_i$ for $i\in [t]$,
where $u_1,\ldots,u_t$ are distinct vertices.

\medskip

\noindent (i) Suppose that $\Delta_t\leq 1$.
If every vertex of degree $1$ is in $N_F[u_1]\cup\cdots\cup N_F[u_{t-1}]$,
then removing $X=\{ u_1,\ldots,u_{t-1}\}$ yields 
three vertices of maximum degree $0$
or a forest of order less than $3$.
Hence, we may assume that $u_t$ is not adjacent to any vertex in $X$.
Now, either removing $X$ yields three vertices of maximum degree $1$,
or removing $X\cup \{ u_t\}$ yields three vertices of maximum degree $0$
or a forest of order less than $3$.
Hence, $\Delta_t\geq 2$.

\medskip

\noindent (ii) Suppose that 
$\Delta_{t+1-i}+2\Delta_{t+2-i}\leq {i+2\choose 2}+2$ for some $i\in [t]\setminus \{ 1\}$.
If $X=\{ u_1,\ldots,u_{t-i}\}$,
then 
$\Delta_1(F-X)+2\Delta_2(F-X)
\leq \Delta_{t+1-i}+2\Delta_{t+2-i}\leq {i+2\choose 2}+2,$ 
and, Theorem \ref{theorem2} implies the contradiction
$f_3(F)\leq (t-i)+f_3(F-X)\leq (t-i)+i=t$,
which completes the proof of (ii).

\medskip

\noindent (iii) By (i) and (ii), we obtain
\begin{eqnarray*} 
\Big(\Delta_1+2\Delta_2\Big)+\Big(\Delta_2+2\Delta_3\Big)+\cdots+\Big(\Delta_{t-1}+2\Delta_t\Big)+\Delta_t
& \geq & \sum\limits_{i=2}^t\left({i+2\choose 2}+3\right)+2\\
& = & \frac{1}{6}t^3 + t^2 + \frac{29}{6}t - 4.
\end{eqnarray*}
Since $\Delta_1\geq 2$, this implies
$3\Big(\Delta_1+\Delta_2+\cdots+\Delta_t\Big)\geq \frac{1}{6}t^3 + t^2 + \frac{29}{6}t$,
which implies (iii).
$\Box$

\medskip

\noindent We obtain a result similar to Theorem \ref{theorem1}.

\begin{corollary}\label{corollary2}
If $t$ is an integer at least $2$, and $F$ is a forest of size less than 
$\frac{1}{18}t^3 + \frac{1}{3}t^2 + \frac{11}{18}t+1$,
then $f_3(F)\leq t$.
\end{corollary}
{\it Proof:} Clearly, we may assume that $F$ has at least $t$ vertices.
Since
$\Delta_1(F)+\Delta_2(F)+\cdots+\Delta_t(F)\leq m(F)+(t-1)$,
Corollary \ref{corollary1}(iii) implies $f_3(F)\leq t$. $\Box$

\medskip

\noindent In order to understand how tight Corollary \ref{corollary2} actually is,
we construct forests $F$ with few edges and a large value of $f_3(F)$.
Therefore, let $a_1=1$, $a_2=3$, and, for every integer $i$ at least $3$, let
\begin{eqnarray}\label{e1}
a_i&=&\max\Big\{ a_{i-1},i-a_{i-1}+2a_{i-2}\Big\}.
\end{eqnarray}
It is easy to verify by induction that $a_{2i+1}=a_{2i}=i^2+i+1$ for every positive integer $i$.

For a positive integer $t$, let $F_t=K_{1,a_1}\cup K_{1,a_2}\cup\cdots\cup K_{1,a_t}$.

\begin{lemma}\label{lemma2}
If $t$ is a positive integer, then 
$f_3(F_t)=t$ and $m(F_t)=\frac{t^3}{12}+O(t^2)$;
more precisely
$$m(F_t)=
\begin{cases}
\frac{2}{3}k^3+2k^2 +\frac{10}{3}k+ 1 & \mbox{, if $t=2k+1$, and}\\
\frac{2}{3}k^3+k^2 +\frac{7}{3}k & \mbox{, if $t=2k$.}
\end{cases}$$
\end{lemma}
{\it Proof:} Since the statement about the size of $F_t$ follows from a straightforward calculation 
using the closed formula for the $a_i$, we only give details for the proof of $f_3(F_t)=t$.
Clearly, removing the $t$ centers of the stars results in an edgeless forest, 
which implies $f_3(F_t)\leq t$.
Now, let $X$ be a minimum set of vertices of $F_t$ such that 
$F_t-X$ has at least three vertices of maximum degree.
Let $\Delta=\Delta_1(F_t-X)$,
and, let $d_{F_t-X}(v_i)=\Delta$ for $i\in [3]$,
where $v_1$, $v_2$, and $v_3$ are distinct.

If $\Delta=0$, then clearly $|X|\geq t$.
Since removing the $t-2$ vertices of largest degree and $2$ endvertices from $K_{1,a_2}$
yields three vertices of maximum degree $1$ in the most efficient way,
if $\Delta=1$, then $|X|\geq (t-2)+2=t$.
Hence, we may assume that $\Delta\geq 2$, 
which implies that $v_1$, $v_2$, and $v_3$ are distinct centers 
of some star components $K_{1,a_i}$ of $F_t$.
Let 
$v_1$ be the center of the component $K_{1,a_p}$,
$v_2$ be the center of the component $K_{1,a_q}$, and
$v_3$ be the center of the component $K_{1,a_r}$, where $p<q<r$.
Clearly, $X$ contains 
$a_r-a_p$ neighbors of $v_3$,
$a_q-a_p$ neighbors of $v_2$,
and at least one vertex from every star component $K_{1,a_i}$ with $q<i<r$ or $r<i\leq t$. Using the monotonicity of the $a_i$ and ($\ref{e1}$), this implies
\begin{eqnarray*}
|X| 
& \geq & (a_r-a_p)+(a_q-a_p)+(r-q-1)+(t-r)\\
& = & a_r+a_q-2a_p+(t-q-1)\\
& \stackrel{mon.}{\geq} & a_{q+1}+a_q-2a_{q-1}+(t-q-1)\\
& \stackrel{(\ref{e1})}{\geq} & (q+1)+(t-q-1)\\
& = & t,
\end{eqnarray*}
which completes the proof. $\Box$
 
\medskip

\noindent Lemma \ref{lemma2} implies that in any version of Corollary \ref{corollary2},
the upper bound on the size is at most $\frac{t^3}{12}+O(t^2)$,
that is, the bound in Corollary \ref{corollary2} might be improved by an asymptotic factor of $3/2$. 
 
The following lemma will be used to extend Theorem \ref{theorem2} to graphs of girth at least $5$ and larger values of $k$.

\begin{lemma}\label{lemma3}
Let $k$ and $t$ be integers with $k\geq 2$ and $t\geq (k-1)^2$.
If $G$ is a graph of girth at least $5$, and 
$$\Delta_1(G)+\cdots+\Delta_{k-1}(G)-(k-1)\Delta_k(G)\leq t,$$
then $f_k(G)\leq t$.
\end{lemma} 
{\it Proof:} Let $\Delta_i=\Delta_i(G)$ and let $d_G(u_i)=\Delta_i$ for $i\in [k]$,
where $u_1,\ldots,u_k$ are distinct vertices.

First, suppose that $\Delta_k<k-1$.
We remove $u_1,\ldots,u_{k-1}$,
and, as long as the current graph has order at least $k$
but less than $k$ vertices of maximum degree,
we iteratively remove all vertices of maximum degree from the current graph.
Therefore,
removing $u_1,\ldots,u_{k-1}$,
at most $(k-1)$ further vertices of degree $k-2$,
at most $(k-1)$ further vertices of degree $k-3$,
and so on, 
until at most $(k-1)$ further vertices of degree $1$,
yields either a graph with $k$ vertices of maximum degree
or a graph with less than $k$ vertices.
Since we removed at most $(k-1)+(k-1)(k-2)=(k-1)^2\leq t$ vertices,
we obtain $f_k(G)\leq t$.
Hence, we may assume that $\Delta_k\geq k-1$.

Let $i\in [k-1]$.
By the girth condition,
$u_i$ has at most $k-1$ neighbors in 
$$N_i=N_G[u_1]\cup \cdots \cup N_G[u_{i-1}]\cup N_G[u_{i+1}]\cup \cdots \cup N_G[u_k].$$
Therefore, there are $\Delta_i-\Delta_k\leq \Delta_i-(k-1)$ neighbors of $u_i$ 
outside of $N_i$ whose removal results in a graph in which $u_i$ has degree $\Delta_k$.
Doing this for every $i$ in $[k-1]$ yields $k$ vertices of maximum degree $\Delta_k$.
$\Box$ 

\medskip

\noindent We proceed to the extension of Theorem \ref{theorem2}.

\begin{theorem}\label{theorem3}
Let $k$ and $t$ be integers with $k\geq 2$ and $t\geq (k-1)^2$.
There is some integer $c_k$ such that, 
if $G$ is a graph of girth at least $5$, and 
$$\Delta_1(G)+2\Delta_2(G)+3\Delta_3(G)+\cdots+(k-1)\Delta_{k-1}(G)\leq {t+2\choose 2}+c_k,$$
then $f_k(G)\leq t$.
\end{theorem}
{\it Proof:} Clearly, we may assume that $G$ has at least $t+k$ vertices.
Let $\Delta_i=\Delta_i(G)$ and let $d_G(u_i)=\Delta_i$ for $i\in [k]$,
where $u_1,\ldots,u_k$ are distinct vertices.
The proof is by induction on $t$.

First, let $t=(k-1)^2$. 
Let $c_k$ be such that ${(k-1)^2+2\choose 2}+c_k=k-1$.
We obtain that $\Delta_1\leq k-1$, and
removing 
at most $(k-1)$ vertices of degree $k-1$,
at most $(k-1)$ further vertices of degree $k-2$,
and so on, 
until at most $(k-1)$ further vertices of degree $1$,
yields either a graph with $k$ vertices of maximum degree
or a graph with less than $k$ vertices.
Since we removed at most $(k-1)^2=t$ vertices,
we obtain $f_k(G)\leq t$.

Next, let $t>(k-1)^2$.
By Lemma \ref{lemma3}, we may assume that 
$\Delta_1(G)+\cdots+\Delta_{k-1}(G)-(k-1)\Delta_k(G)\geq t+1$.
Similarly as in the proof of Theorem \ref{theorem2},
we may assume, by induction, that 
$$\Delta_2(G)+2\Delta_3(G)+3\Delta_4(G)+\cdots+(k-1)\Delta_k(G)\geq {t-1+2\choose 2}+c_k+1.$$
Adding these two inequalities implies a contradiction,
which completes the proof. $\Box$

\medskip

\noindent Theorem \ref{theorem3} has several interesting consequences.

\begin{corollary}\label{corollary3}
Let $k$ be a fixed integer at least $2$.

There is a function $g:\mathbb{N}\to \mathbb{Z}$ with $|g(t)|=O(t^2)$ such that, 
if $t$ is some positive integer, and 
$G$ is a graph of size at most $\frac{t^3}{6{k\choose 2}}+g(t)$ and girth at least $5$,
then $f_k(G)\leq t$.
\end{corollary}
{\it Proof:} Choosing $g(t)$ equal to $-\frac{t^3}{6{k\choose 2}}$ for $t<(k-1)^2$,
the statement becomes trivial for $t<(k-1)^2$.
Hence, we may assume that $t\geq (k-1)^2$.

Let the graph $G$ of girth at least $5$ be such that $f_k(G)>t$; 
in particular, $G$ has at least $t+k$ vertices.
Let $\Delta_i=\Delta_i(G)$ for $i\in [t]$.
Arguing similarly as in the proof of Corollary \ref{corollary1} (ii), we obtain that
$$\Delta_{t+1-i}+2\Delta_{t+2-i}+\cdots+(k-1)\Delta_{t+k-1-i}\geq {i+2\choose 2}+c_k+1$$
for every $i\in [t]\setminus \Big[(k-1)^2-1\Big]$.
Adding all these inequalities, we obtain, using $1+2+\cdots+(k-1)={k\choose 2}$, that
\begin{eqnarray*}
{k\choose 2}\Big(\Delta_1+\cdots+\Delta_{t+k-1-(k-1)^2}\Big) & \geq & 
\sum\limits_{i=(k-1)^2}^t\left({i+2\choose 2}+c_k+1\right)
=\frac{t^3}{6}+O(t^2),
\end{eqnarray*}
where the implicit constants depend on the fixed value of $k$.

If $H$ is the subgraph of $G$ induced by the $t+k-1-(k-1)^2<t$ vertices of the largest degrees,
then
$$m(G)\geq 
\Big(\Delta_1+\cdots+\Delta_{t+k-1-(k-1)^2}\Big)-m(H)
\geq \frac{t^3}{6{k\choose 2}}+O(t^2),$$
which completes the proof. $\Box$

\medskip

\noindent It is a simple consequence of the Moore bound \cite{ms} that,
for every positive integer $p$,
we have $m(G)\leq 2n(G)^{\frac{p+1}{p}}$
for every graph $G$ of girth more than $2p$.

\begin{corollary}\label{corollary4}
Let $k$ and $p$ be fixed integers with $k\geq 2$ and $p\geq 3$.

If $G$ has girth more than $2p$, then
$$f_k(G)\leq \Big(1+o(1)\Big)\left(12{k\choose 2}\right)^{\frac{1}{3}}n(G)^{\frac{p+1}{3p}}.$$
\end{corollary}
{\it Proof:} Let $G$ be a graph of girth more than $2p$, and let $t=f_k(G)-1$.
By the above consequence of the Moore bound and Corollary \ref{corollary3}, we obtain
\begin{eqnarray*}
n(G) & \geq & \left(\frac{1}{2}m(G)\right)^{\frac{p}{p+1}}
> \left(\frac{1}{2}\left(\frac{1}{6{k\choose 2}}+o(1)\right)t^3\right)^{\frac{p}{p+1}}
= \left(\left(\frac{1}{12{k\choose 2}}+o(1)\right)t^3\right)^{\frac{p}{p+1}}.
\end{eqnarray*}
This implies
$t<\Big(1+o(1)\Big)\Big(12{k\choose 2}\Big)^{\frac{1}{3}}n(G)^{\frac{p+1}{3p}}$,
which completes the proof. $\Box$

\medskip

\noindent Arguing in a similar way for forests, we obtain the following considerable improvement of Theorem \ref{theoremclz}.

\begin{corollary}\label{corollary5}
Let $k$ be a fixed integer with $k\geq 2$.

If $F$ is a forest, then
$$f_k(G)\leq \Big(1+o(1)\Big)\left(6{k\choose 2}\right)^{\frac{1}{3}}n(G)^{\frac{1}{3}}.$$
\end{corollary}

\section{An algorithm for forests}

In this section we describe an efficient algorithm calculating $f_k(F)$ for a given forest $F$.

Let $k$ be an integer at least $2$.
Let $T$ be a tree of order more than $k$, 
let $S$ be a set of $k$ distinct vertices of $T$, and, 
let $\Delta$ be some non-negative integer at most $\Delta(T)$.
The vertices in $S$ are called {\it special}.
We root $T$ in some non-special vertex $r$,
and, for every vertex $u$ of $T$, 
we denote by $T(u)$ the subtree of $T$ rooted in $u$
and containing $u$ as well as all descendants of $u$.

For a vertex $u$ of $T$, let $(n_1(u),n_2(u),n_3(u))$ be a triple of integers, where
\begin{enumerate}[(i)]
\item $n_1(u)$ is the maximum order of an induced subforest $T_1(u)$ of $T(u)$ such that 
\begin{itemize}
\item $u\not\in V(T_1(u))$, 
\item $S\cap V(T(u))\subseteq V(T_1(u))$,
\item $\Delta(T_1(u))\leq \Delta$, and
\item $d_{T_1(u)}(v)=\Delta$ for every vertex $v\in S\cap V(T(u))$.
\end{itemize}
Note that, if $u$ is special, then $n_1(u)=\max\emptyset$, which, by convention, is $-\infty$.
\item $n_2(u)$ is the maximum order of an induced subforest $T_2(u)$ of $T(u)$ such that 
\begin{itemize}
\item $\{ u\} \cup \Big(S\cap V(T(u))\Big)\subseteq V(T_2(u))$,
\item $\Delta(T_2(u))\leq \Delta$, and
\item $d_{T_2(u)}(v)=\Delta$ for every vertex $v\in \{ u\}\cup \Big(S\cap V(T(u))\Big)$.
\end{itemize}
\item $n_3(u)$ is the maximum order of an induced subforest $T_3(u)$ of $T(u)$ such that 
\begin{itemize}
\item $\{ u\} \cup \Big(S\cap V(T(u))\Big)\subseteq V(T_3(u))$,
\item $\Delta(T_3(u))\leq \Delta$, 
\item $d_{T_3(u)}(v)=\Delta$ for every vertex $\Big(S\cap V(T(u))\Big)\setminus \{ u\}$, and
\item if $u$ is special, then $d_{T_3(u)}(u)=\Delta-1$,
and, if $u$ is non-special, then $d_{T_3(u)}(u)\leq \Delta-1$.
\end{itemize}
\end{enumerate}
If $u$ is a non-special leaf of $T$, then 
$$(n_1(u),n_2(u),n_3(u))=
\begin{cases}
(0,1,-\infty) & \mbox{, if $\Delta=0$ and}\\
(0,-\infty,1) & \mbox{, if $\Delta\geq 1$,}
\end{cases}$$
and, if $u$ is a special leaf of $T$, then 
$$(n_1(u),n_2(u),n_3(u))=
\begin{cases}
(-\infty,1,-\infty) & \mbox{, if $\Delta=0$,}\\
(-\infty,-\infty,1) & \mbox{, if $\Delta=1$, and}\\
(-\infty,-\infty,-\infty) & \mbox{, if $\Delta\geq 2$.}
\end{cases}$$
The following lemma gives recursions for non-leaf vertices of $T$. 

\begin{lemma}\label{lemma4}
Let $u$ be a non-leaf vertex of $T$, where we use the notation introduced above.

Let $v_1,\ldots,v_p$ be the special children of $u$,
and, let $w_1,\ldots,w_q$ be the non-special children of $u$.

Let 
$$n_3(w_1)-n_1(w_1)\geq n_3(w_2)-n_1(w_2)\geq \ldots \geq n_3(w_q)-n_1(w_q).$$
If $n_3(w_1)-n_1(w_1)<0$, let $q'=0$, and, 
if $n_3(w_1)-n_1(w_1)\geq 0$, 
let $q'\in [q]$ be maximum such that $n_3(w_{q'})-n_1(w_{q'})\geq 0$.
\begin{enumerate}[(i)]
\item If $u$ is non-special, then 
$$n_1(u)=
\sum_{i=1}^p n_2(v_i)+\sum_{j=1}^q\max\Big\{ n_1(w_j),n_2(w_j),n_3(w_j)\Big\}.$$
\item If $p>\Delta$ or $p+q<\Delta$, then $n_2(u)=-\infty$, and, if $p\leq \Delta\leq p+q$, then
$$n_2(u)=
\sum_{i=1}^p n_3(v_i)
+\sum_{j=1}^{\Delta-p}n_3(w_j)
+\sum_{j=\Delta-p+1}^{q}n_1(w_j).$$
\item If $u$ is special, and $p>\Delta-1$ or $p+q<\Delta-1$, then $n_3(u)=-\infty$, 
and, if $u$ is special, and $p\leq \Delta-1\leq p+q$, then  
$$n_3(u)=
\sum_{i=1}^p n_3(v_i)
+\sum_{j=1}^{\Delta-1-p}n_3(w_j)
+\sum_{j=\Delta-p}^{q}n_1(w_j).$$
\item If $u$ is non-special and $p>\Delta-1$, then $n_3(u)=-\infty$, 
and, if $u$ is non-special and $p\leq \Delta-1$, then  
$$n_3(u)=
\sum_{i=1}^p n_3(v_i)
+\sum_{j=1}^{\min\{ q',\Delta-1-p\}}n_3(w_j)
+\sum_{j=\min\{ q',\Delta-1-p\}+1}^{q}n_1(w_j).$$
\end{enumerate}
\end{lemma}
{\it Proof:} (i) Since $u$ does not belong to $T_1(u)$,
for every special child $v_i$ of $u$, 
the forest $T_1(u)\cap T(v_i)$
has at most as many vertices as $T_2(v_i)$,
and,
for every non-special child $w_j$ of $u$, 
the forest $T_1(u)\cap T(w_j)$
has at most as many vertices as 
the forest of largest order in $\{ T_1(w_j),T_2(w_j),T_3(w_j)\}$,
which implies that $n_1(u)$ is at most the specified value.
On the other hand, 
combining the mentioned forests in the obvious way, 
it follows that $n_1(u)$ is also at least the specified value,
which completes the proof of (i).

\medskip

\noindent (ii) If $p>\Delta$ or $p+q<\Delta$, 
then no forest with the properties required for $T_2(u)$
exists, and, hence, $n_2(u)=-\infty$.
If $p\leq \Delta\leq p+q$,
then, since $u$ belongs to $T_2(u)$ and has degree exactly $\Delta$ in $T_2(u)$,
for every special child $v_i$ of $u$, 
the forest $T_2(u)\cap T(v_i)$
has at most as many vertices as $T_3(v_i)$,
there are exactly $\Delta-p$ non-special children $w_j$ of $u$ 
that belong to $T_2(u)$,
and the forest $T_2(u)\cap T(w_j)$ for such a $w_j$
has at most as many vertices as $T_3(w_j)$,
and, for the remaining $q-(\Delta-p)$ non-special children $w_j$ of $u$ 
that do not belong to $T_2(u)$,
the forest $T_2(u)\cap T(w_j)$ for such a $w_j$
has at most as many vertices as $T_1(w_j)$.
In view of the ordering of the non-special children $w_j$ of $u$,
this implies that $n_2(u)$ has at most the specified value.
Again, on the other hand, 
combining the mentioned forests in the obvious way, 
it follows that $n_2(u)$ is also at least the specified value, 
which completes the proof of (ii).

\medskip

\noindent Since the proof of (iii) is almost identical to the proof of (ii), 
we proceed to the proof of (iv).
Since $u$ belongs to $T_3(u)$ and has degree at most $\Delta-1$,
some non-special child $w_j$ of $u$ may only belong to $T_3(u)$ 
if $n_3(w_j)-n_1(w_j)\geq 0$, 
which easily implies (iv) arguing similarly as for the proof of (ii). $\Box$

\begin{theorem}\label{theorem4}
For a fixed integer $k$ at least $2$, and a given forest $F$, the value $f_k(F)$
can be determined in polynomial time.
\end{theorem}
{\it Proof:} Clearly, we may assume that $F$ has more than $k$ vertices.
Let $S$ be a set of $k$ distinct vertices of $F$, and let $\Delta$ be a non-negative integer at most $\Delta(F)$.
If $F$ is disconnected, then we add a vertex $r$ with a neighbor in each component of $F$,
and denote the resulting tree by $T$.
Otherwise, let $T=F$, and let $r$ be a vertex of $F$ that does not belong to $S$.
Using the recursions from Lemma \ref{lemma4},
we can determine, in polynomial time, 
$(n_1(r),n_2(r), n_3(r))$ for $T$,
denoted by 
$\left(n_1^{(S,\Delta)}(r),n_2^{(S,\Delta)}(r),n_3^{(S,\Delta)}(r)\right)$
for this specific choice of $S$ and $\Delta$.

If $F$ is connected, then $n(F)-f_k(F)$ equals
$$
\max\left\{
\max\left\{ n_1^{(S,\Delta)}(r),n_2^{(S,\Delta)}(r),n_3^{(S,\Delta)}(r)\right\}:
S\in{ V(F)\choose k}\mbox{ and }\Delta\in \{ 0\}\cup [\Delta(F)]\right\},$$
and,
if $F$ is not connected, then $n(F)-f_k(F)$ equals
$$
\max\left\{ n_1^{(S,\Delta)}(r):
S\in{ V(F)\choose k}\mbox{ and }\Delta\in \{ 0\}\cup [\Delta(F)]\right\}.$$
Since these maxima are taken over polynomially many values, 
the desired statement follows.
$\Box$

\medskip

\noindent It seems possible yet challenging to extend this approach to graphs of bounded tree width.


\begin{thebibliography}{99}
\bibitem{ab} M.O. Albertson, D.L. Boutin, Lower bounds for constant degree independent sets, Discrete Mathematics 127 (1994) 15-21.
\bibitem{bs} B. Bollob\'{a}s, A.D. Scott, Independent sets and repeated degrees, Discrete Mathematics 170 (1997) 41-49.
\bibitem{clz} Y. Caro, J. Lauri, C. Zarb, Equating two maximum degrees, arXiv:1704.08472v1.
\bibitem{csy} Y. Caro,  A. Shapira, R. Yuster, Forcing $k$-repetitions in degree sequences,
The Electronic Journal of Combinatorics 21(2014) P1-24.
\bibitem{cy} Y. Caro, R. Yuster, Large induced subgraphs with equated maximum degree, Discrete Mathematics 310 (2010) 742-747.
\bibitem{cw} Y. Caro, D.B. West, Repetition number of graphs, The Electronic Journal of Combinatorics 16 (2009) R7.
\bibitem{ms} M. Miller, J. \v{S}ir\'{a}\v{n}, Moore graphs and beyond: A survey of the degree/diameter problem, The Electronic Journal of Combinatorics 20 (2013) $\#$ DS14v2.
\end{thebibliography}
\end{document}